\documentclass[a4paper,12pt]{amsart}
\usepackage{amssymb}
\usepackage{ifthen}
\usepackage{cite}
\usepackage[dvips]{graphicx}
\nonstopmode \numberwithin{equation}{section}
\usepackage{amssymb}
\usepackage{ifthen}
\usepackage{graphicx}
\usepackage{amsmath}
\usepackage[T1]{fontenc} 
\usepackage[utf8]{inputenc}
\usepackage[usenames,dvipsnames]{color}
\usepackage{color}
\usepackage[english]{babel}
\usepackage{fancyhdr}
\usepackage{fancybox}
\usepackage{tikz}

\nonstopmode \numberwithin{equation}{section}
\theoremstyle{plain}

\newtheorem{conj}{Conjecture}

\theoremstyle{definition}
\newtheorem{defn}{Definition}[section]

\newtheorem{thm}{Theorem}[section]
\newtheorem{prob}{Problem}[section]
\newtheorem{cor}{Corollary}[section]
\newtheorem{ques}{Question}[section]
\newtheorem{prop}{Proposition}[section]
\newtheorem{rem}{Remark}[section]
\newtheorem{lem}{Lemma}[section]

\newcounter{minutes}
\divide\time by 60
\newcounter{hours}
\multiply\time by 60
\addtocounter{minutes}{-\time}

\newcounter {own}
\def\theown {\thesection       .\arabic{own}}

\newenvironment{pf}[1][]{%
 \vskip 3mm
 \noindent
 \ifthenelse{\equal{#1}{}}%
  {{\slshape Proof. }}%
  {{\slshape #1.} }%
 }%
{\qed\bigskip}

\newcounter{alphabet}

\def\be{\begin{equation}}
\def\ee{\end{equation}}

\newcommand{\bee}{\begin{enumerate}}
\newcommand{\eee}{\end{enumerate}}

\newcommand{\blem}{\begin{lem}}
\newcommand{\elem}{\end{lem}}
\newcommand{\bthm}{\begin{thm}}
\newcommand{\ethm}{\end{thm}}
\newcommand{\bcor}{\begin{cor}}
\newcommand{\ecor}{\end{cor}}
\newcommand{\beg}{\begin{examp}}
\newcommand{\eeg}{\end{examp}}
\newcommand{\begs}{\begin{examples}}
\newcommand{\eegs}{\end{examples}}

\newcommand{\bdefn}{\begin{defn}}
\newcommand{\edefn}{\end{defn}}

\newcommand{\bprob}{\begin{prob}}
\newcommand{\eprob}{\end{prob}}
\newcommand{\bei}{\begin{itemize}}
\newcommand{\eei}{\end{itemize}}

\newcommand{\bcon}{\begin{conj}}
\newcommand{\econ}{\end{conj}}
\newcommand{\bcons}{\begin{conjs}}
\newcommand{\econs}{\end{conjs}}
\newcommand{\bprop}{\begin{prop}}
\newcommand{\eprop}{\end{prop}}
\newcommand{\br}{\begin{rem}}
\newcommand{\er}{\end{rem}}
\newcommand{\brs}{\begin{rems}}
\newcommand{\ers}{\end{rems}}
\newcommand{\bo}{\begin{obser}}
\newcommand{\eo}{\end{obser}}
\newcommand{\bos}{\begin{obsers}}
\newcommand{\eos}{\end{obsers}}
\newcommand{\bpf}{\begin{pf}}
\newcommand{\epf}{\end{pf}}
\newcommand{\ba}{\begin{array}}
\newcommand{\ea}{\end{array}}
\newcommand{\beq}{\begin{eqnarray}}
\newcommand{\beqq}{\begin{eqnarray*}}
\newcommand{\eeq}{\end{eqnarray}}
\newcommand{\eeqq}{\end{eqnarray*}}

\begin{document}
\title{Generalized Bohr inequalities for certain classes of functions and their applications}

\author{Molla Basir Ahamed}
\address{Molla Basir Ahamed, Department of Mathematics, Jadavpur University, Kolkata-700032, West Bengal,India.}
\email{mbahamed.math@jadavpuruniversity.in}

\subjclass[{AMS} Subject Classification:]{Primary 47A56; 30B10; 47A63; 30C80, 30A10, 30C35, 30C62, 31A05, Secondary 30C45, 30H05}

\def\thefootnote{}
\footnotetext{ {\tiny File:~\jobname.tex,
printed: \number\year-\number\month-\number\day,
          \thehours.\ifnum\theminutes<10{0}\fi\theminutes }
} \makeatletter\def\thefootnote{\@arabic\c@footnote}\makeatother

\begin{abstract} 
Let $ \mathcal{B}:=\{f(z)=\sum_{n=0}^{\infty}a_nz^n\; \mbox{with}\; |f(z)|<1\;\mbox{for all}\; z\in\mathbb{D}\} $. The improved version of the classical Bohr's inequality \cite{Bohr-1914} states that if $ f\in\mathcal{B} $, then the associated majorant series $ M_f(r):=\sum_{n=0}^{\infty}|a_n|r^n\leq 1 $ holds for $ |z|=r\leq 1/3 $ and the constant $ 1/3 $ cannot be improved. Bohr's original theorem and its subsequent generalizations remain active fields of study, driving investigations in a wide range of function spaces. In this paper, first we establish a generalized Bohr inequality for the class $\mathcal{B}$ by allowing a sequence $\{\varphi_n(r)\}_{n=0}^{\infty}$ of non-negative continuous functions on $[0, 1)$ in the place of $\{r^n\}_{n=0}^{\infty}$ of the majorant series $M_f(r)$ introducing a weighted sequence of non-negative continuous functions $\{\Phi_n(r)\}_{n=0}^{\infty}$ on $[0, 1)$. Secondly, as a generalization, we obtain a refined version of the Bohr inequality for a certain class $\tilde{G}^0_{\mathcal{H}}(\beta)$ of harmonic mappings. All the results are proved to be sharp.
\end{abstract}

\maketitle
\pagestyle{myheadings}
\markboth{M. B. Ahamed}{Generalized Bohr inequalities for certain classes of functions and their applications}
\tableofcontents
\section{Introduction}
Bohr's classical theorem and its generalizations are now active areas
of research and have been the source of investigations in numerous function spaces. Let $ \mathcal{A} $ denote the space of all functions analytic in the unit disk $ \mathbb{D}:=\{z\in\mathbb{C} : |z|<1\} $
equipped with the topology of uniform convergence on compact subsets of $ \mathbb{D} $. We define $ \mathcal{B}:=\{f\in\mathcal{A} : f(z)=\sum_{n=0}^{\infty}a_nz^n\; \mbox{with}\; |f(z)|<1\;\mbox{for all}\; z\in\mathbb{D}\} $. The improved version of the classical Bohr's inequality \cite{Bohr-1914} states that if $ f\in\mathcal{B} $, then the associated majorant series $ M_f(r):=\sum_{n=0}^{\infty}|a_n|r^n\leq 1 $ holds for $ |z|=r\leq 1/3 $ and the constant $ 1/3 $ cannot be improved. The constant $ 1/3 $ is famously known as the Bohr radius and the inequality $ M_f(r)\leq 1 $ is known as the Bohr inequality for the class $ \mathcal{B} $. Henceforth, for a class $ \mathcal{F} $, if there exists a positive real number $ r_0 $ such that $ M_f(r)\leq 1 $ for $ |z|=r\leq r_0 $ and for $ f\in\mathcal{F} $, then we say that $ r_0 $ is the Bohr radius for the class $ \mathcal{F} $. Moreover, the constant $ r_0 $ is said to be best possible in the sense that for a function $ f_0\in\mathcal{F} $, if $ M_{f_0}(r)>1 $ whenever $ r>r_0 $. \vspace{1.2mm}

In the literature, Bohr’s power series theorem has been studied in many different situations which are called the Bohr phenomenon. In the study of the Bohr phenomenon, the initial term $ |a_0| $ in $ M_f(r) $ plays some crucial role. For instance, Bohr's theorem has been studied by introducing the term $ |a_0|^p $, where $ 0<p\leq 2 $, instead of $ |a_0| $, with the corresponding radius $ p/(p+2) $ (see \cite{Liu-Ponnusamy-PAMS-2021} and references therein). In addition, if  $|a_0|$ is replaced by $|f(z)|$, then the constant $1/3$ could be replaced by $\sqrt{5}-2$ which is best possible (see  \cite{Alkhaleefah-Kayumov-Ponnusamy-PAMS-2019}). A detailed account of research on the Bohr radius problem can be found in the survey article \cite{Ponnusmy-Survey} and references therein.\vspace{1.2mm}

 Actually, Bohr’s theorem received greater interest after it was used by Dixon \cite{Dixon & BLMS & 1995} to characterize Banach algebras that satisfy von Neumann's inequality. The generalization of Bohr’s theorem is now an active area of research: for instance, Aizenberg \textit{et al.} \cite{aizenberg-2001}, and Aytuna and Djakov \cite{Aytuna-Djakov-BLMS-2013} have studied the Bohr property of bases for holomorphic functions; Ali \textit{et al.} \cite{Ali-Abdul-NG-CVEE-2016} have found the Bohr radius for the class of starlike log-harmonic mappings; while Paulsen \textit{et al.} \cite{Paulsen-PLMS-2002} extended the Bohr inequality to Banach algebras; Hamada \emph{et al.}\cite{Hamada-IJM-2009} have studied the Bohr's theorem for holomorphic mappings with values in homogeneous balls; Galicer \emph{et al.}\cite{Galicer-Mansilla-Muro-TAMS-2020} have studied mixed Bohr radius in several complex variables. For recent development in Bohr inequalities for different class of functions, we refer to the articles  \cite{Aha-Aha-CMFT-2023,Aha-Allu-RMJ-2022,Alkhaleefah-Kayumov-Ponnusamy-PAMS-2019,Allu-Arora-JMAA-2022,Allu-CMB-2022,Das-JMAA-2022,Lata-Singh-PAMS-2022,S. Kumar-PAMS-2022,Kumar-JMAA-2023,Aizen-PAMS-2000,Evdoridis-Ponn-RM-2021,Ponnusamy-Vij-Wirth-JMAA-2022,Aizn-ST-2007,Kumar-CVEE-2023,Kayu-Kham-Ponn-MJM-2021,Boas-1997,Wu-Wang-Long-RACSAM-2022,Gang-Liu-JMAA-2021} and references therein. Study of Bohr inequality also get much attention because of its generalization to holomorphic function in several complex variables established by Boas and Khavinson \cite{Boas-1997} and this result generates an extensive research activity what is called multidimensional Bohr phenomenon. For recent developments in multidimensional Bohr, we refer to articles \cite{S. Kumar-PAMS-2022,Lin-Liu-Ponnusamy-Acta-2023,Kumar-JMAA-2023,Liu-Ponnusamy-PAMS-2021} and references therein. A logarithmic lower bound for multi-dimensional Bohr radii is established in \cite{Defant-Frerick-IJM-2001}, and subsequently, an improved lower bound for the multidimensional Bohr radius is obtained in \cite{Das-CMB-2023}.\vspace{1.2mm}

 However, it can be noted that not every class of functions has the Bohr phenomenon, for example, B\'an\'at\'aau \emph{et al.} \cite{Beneteau-2004} showed that there is no Bohr phenomenon in the Hardy space $ H^p(\mathbb{D},X), $ where $p\in [1,\infty).$ In \cite{Liu-Liu-JMAA-2020}, Liu and Liu have shown that Bohr's inequality fails to hold for the class $ \mathcal{H}(\mathbb{D}^2, \mathbb{D}^2) $, a set of holomorphic functions $ f : \mathbb{D}^2\rightarrow \mathbb{D}^2 $ having lacunary series expansion. The study of the Bohr phenomenon in relation to operator-valued functions has gained attention in recent years, resulting in significant findings (see \cite{Bhowmik-Das-PEMS-2021,Allu-Hal-CMB-2022,Allu-Hal-CMB-2023} and references therein).\vspace{1.2mm}
 
 The objective of this paper is to study generalized Bohr inequalities and their applications for two classes of functions: analytic functions and harmonic mappings. We will analyze their generalizations in relation to existing Bohr inequalities, encompassing both refined and improved variants.
\section{Generalization of Bohr inequality with the sequence of non-negative continuous functions $ \{\varphi_k(r)\}_{k=0}^{\infty} $}
The generalization of the Bohr radius is an important study in understanding the Bohr phenomenon for various classes of functions. The idea involves replacing the sequence $ \{r^k\}_{k=0}^{\infty} $ by a more general setting in the majorant series $M_f(r)$ of the classical Bohr inequality for the class $ \mathcal{B} $ with a generalized sequence. Let $ \{\varphi_k(r)\}_{k=0}^{\infty} $ be a sequence of non-negative continuous functions in $ [0, 1) $ such that the series $ \sum_{k=0}^{\infty}\varphi_k(r) $ converges locally uniformly with respect to $ r\in [0, 1) $. The Bohr radius generalizes in \cite{Kayu-Kham-Ponn-MJM-2021} with the setting og $ \{\varphi_k(r)\}_{k=0}^{\infty} $. Moreover, Chen \emph{et al.} \cite{chen-Liu-Ponnusamy-arXiv-2023} established several significant result corresponding to Bohr-type inequalities for unimodular bounded analytic functions in the unit disk by allowing the sequence $\{\varphi_k(r)\}_{k=0}^{\infty}$ in place of the $\{r^k\}_{k=0}^{\infty}$ in the power series representations of the functions involved with the Bohr sums. Upon an in-depth study of the results in \cite{Kayu-Kham-Ponn-MJM-2021}, however, we observe that the generalizations are significant but possess a limited scope of applications. To elaborate further, it should be specified that the applications derived from those findings do not account for the improved Bohr inequalities (see \cite{Kayumov-CRACAD-2018}) and the refined version of the Bohr inequalities (see \cite{Ponnusamy-Vijayak-Wirths-RM-2020,Liu-Liu-Ponnusamy-2021,Kayumov-Ponnusamy-JMAA-2018}).\vspace{1.2mm}

We study the Bohr inequality with a weighted sequence $\{\Phi_k(r)\}_{k=0}^{\infty}$ of non-negative continuous functions on $[0, 1)$ associated with the coefficients in the quantity
\begin{align*}
	||f_0||^2_r=\sum_{k=1}^{\infty}|a_k|^2r^{2k},\; \mbox{where}\; f_0(z)=f(z)-f(0).
\end{align*}
The following is the main result of this section. The subsequent applications of this result will encompass both the improved and refined Bohr inequalities.
\begin{thm}\label{B-thm-4.2}
	Let $ \{\varphi_k(r)\}_{k=0}^{\infty} $ and $ \{\Phi_k(r)\}_{k=0}^{\infty} $ be sequences of non-negative continuous functions in $ [0, 1) $ such that the series 
	\begin{align*}
		\varphi_0(r)+\sum_{k=1}^{\infty}\varphi_k(r)\; \mbox{and}\; \sum_{k=1}^{\infty}\Phi_k(r)\{\varphi_k(r)\}^2
	\end{align*} 
	both converge locally uniformly with respect to $ r\in [0, 1) $. Let $ f\in\mathcal{B} $ have the series representation $ f(z)=\sum_{k=0}^{\infty}a_kz^k $ in $ \mathbb{D} $ and $ p\in (0, 2\pi] $. If 
	\begin{align*}
		\varphi_0(r)>\frac{2}{p}\sum_{k=1}^{\infty}\left(\varphi_k(r)+\Phi_k(r)\left(\varphi_k(r)\right)^2\right),
	\end{align*}
	where $ R $ is the minimal root of the equation  
	\begin{align*}
		\varphi_0(r)=\frac{2}{p}\sum_{k=1}^{\infty}\left(\varphi_k(r)+\Phi_k(r)\left(\varphi_k(r)\right)^2\right),
	\end{align*} then the following sharp inequality holds:
	\begin{align*}
		\mathcal{C}_f(\varphi, \Phi_k, p, r):=|a_0|^p\varphi_0(r)+\sum_{k=1}^{\infty}|a_k|\varphi_k(r)+\sum_{k=1}^{\infty}\Phi_k(r)|a_k|^2\left(\varphi_k(r)\right)^2\leq \varphi_0(r)
	\end{align*}
	for all $ |z|=r\leq R $, where $ R $ is the minimal positive root of equation
	\begin{align*}
		\varphi_0(r)=\frac{2}{p}\sum_{k=1}^{\infty}\left(\varphi_k(r)+\Phi_k(r)\left(\varphi_k(r)\right)^2\right)
	\end{align*}
	In the case when 
	\begin{align*}
		\varphi_0(r)<\frac{2}{p}\sum_{k=1}^{\infty}\left(\varphi_k(r)+\Phi_k(r)\left(\varphi_k(r)\right)^2\right)
	\end{align*}
	in some interval $ (R, R+\epsilon) $, the number $ R $ cannot be improved. If the the function $ \varphi_k(r) $ $ (k\geq 0) $ are smooth function then the last condition is equivalent to the inequality 
	\begin{align*}
		\varphi^{\prime}_0(R)<(2/p)\sum_{k=1}^{\infty}\left(\varphi^{\prime}_k(R)+\left(\Phi_k(R)\left(\varphi_k(R)\right)^2\right)^{\prime}\right).
	\end{align*}
\end{thm}
As a consequence of Theorem \ref{B-thm-4.2}, allowing $\Phi_k(r)=0$ for each $k\in\mathbb{N}$, we see that the following result can be obtained.
\begin{cor}\cite[Theorem 1]{Kayu-Kham-Ponn-MJM-2021}\label{Cor-2.1}
 Let $ f\in\mathcal{B} $, $ f(z)=\sum_{k=0}^{\infty}a_kz^k $ and $ p\in (0, 2\pi] $. If 
\begin{align*}
\varphi_0(r)>\frac{2}{p}\sum_{k=1}^{\infty}\varphi_k(r)\; \mbox{for}\; r\in [0, R_1),
\end{align*}
where $ R_1 $ is the minimal root of the equation  $ \varphi_0(r)=(2/p)\sum_{k=1}^{\infty}\varphi_k(r) $, then the following sharp inequality holds:
\begin{align*}
A_f(\varphi, p, r)=|a_0|^p\varphi_0(r)+\sum_{k=1}^{\infty}|a_k|\varphi_k(r)\leq \varphi_0(r)\; \mbox{for all}\; r\leq R_1.
\end{align*}
In the case when $ \varphi_0(r)<(2/p)\sum_{k=1}^{\infty}\varphi_k(r) $ in some interval $ (R_1, R_1+\epsilon) $, the number $ R_1 $ cannot be improved. If the the function $ \varphi_k(r) $ $ (k\geq 0) $ are smooth function then the last condition is equivalent to the inequality $ \varphi^{\prime}_0(R_1)<(2/p)\sum_{k=1}^{\infty}\varphi^{\prime}_k(R_1) $. 
\end{cor}
\begin{rem} It can be readily seen from the following observations:
	\begin{enumerate}
\item[(i)] In particular, if $ \Phi_k(r)=0 $ on $(0, 1)$, then it is easy to see that Corollary \ref{B-thm-4.2} reduces to exactly Corollary \ref{Cor-2.1}. As a consequence, Theorem \ref{B-thm-4.2} provides a generalization of \cite[Theorem 1]{Kayu-Kham-Ponn-MJM-2021}, facilitating the ease of implementing the corresponding applications outlined in \cite{Kayu-Kham-Ponn-MJM-2021}.
\item[(ii)] Moreover, if $\Phi_k(r) > 0$ on $(0, 1)$, then Theorem \ref{B-thm-4.2} becomes an improved and refined version of \cite[Theorem 1]{Kayu-Kham-Ponn-MJM-2021}, which will be discussed below as applications.
	\end{enumerate}
\end{rem}
\subsection{Some applications of Theorem \ref{B-thm-4.2}}
Suppose that $f\in\mathcal{B}$, $f(z)=\sum_{k=0}^{\infty}a_kz^k$ and $p\in [0, 2)$. Then Theorem \ref{B-thm-4.2} yields the following:
	\begin{enumerate}
	\item[(I)] {\bf Refined Bohr inequalities}:

		\item[(i)] For $\varphi_k(r)=r^k$ $(k\geq 0)$ and $\Phi_k(r)=\left(\frac{1}{1+a}+\frac{r}{1-r}\right)$, $(k\geq 1)$, where $a=|a_0|=|f(0)|$, we obtain (see \cite[Theorem 2]{Ponnusamy-Vijayak-Wirths-RM-2020}):
		\begin{enumerate}
			\item[(a)] if $p=1$, then
			\begin{align*}
				\sum_{k=0}^{\infty}|a_k|r^k+\left(\frac{1}{1+a}+\frac{r}{1-r}\right)\sum_{k=1}^{\infty}|a_k|^2r^{2k}\leq 1
			\end{align*}
			for $r\leq {1}/{(2+a)}$. The numbers ${1}/{(2+a)}$ and ${1}/{(1+a)}$ cannot be improved.
			\item[(b)] if $p=2$, then
			\begin{align*}
				a^2+\sum_{k=1}^{\infty}|a_k|r^k+\left(\frac{1}{1+a}+\frac{r}{1-r}\right)\sum_{k=1}^{\infty}|a_k|^2r^{2k}\leq 1
			\end{align*}
			for $r\leq {1}/{2}$ and the numbers ${1}/{2}$ and ${1}/{(1+a)}$ cannot be improved.
		\end{enumerate}
		\item[(ii)] For 
		\begin{align*}
			\varphi_k(r)=
			\begin{cases}
				0,\; \mbox{for}\; k=0\\
				r^k,\; \mbox{for}\; k\geq 1
			\end{cases}\; \mbox{and}\;\Phi_k(r)=
			\begin{cases}
				\dfrac{-1}{|a_1|^2r^2},\; \mbox{for}\; k=1\vspace{2mm}\\
				\left(\dfrac{r^{-1}}{1+|a_1|}+\dfrac{1}{1-r}\right),\; \mbox{for}\; k\geq 2,
			\end{cases}
		\end{align*}
		we obtain (see \cite[Theorem 3(a)]{Ponnusamy-Vijayak-Wirths-RM-2020}):
		\begin{align*}
			\sum_{k=1}^{\infty}|a_k|r^k+\left(\frac{1}{1+|a_1|}+\frac{r}{1-r}\right)\sum_{k=1}^{\infty}|a_k|^2r^{2k-1}\leq 1
		\end{align*}
		for $r\leq {3}/{5}$ and the number ${3}/{5}$ is sharp.
		\item[(iii)]
		For 
		\begin{align*}
			\varphi_k(r)=
			\begin{cases}
				0,\; \mbox{for}\; k=0\\
				r^k,\; \mbox{for}\; k\geq 1
			\end{cases}
		\end{align*}
		and 
		\begin{align*}
			\Phi_k(r)=
			\begin{cases}
				\left(\dfrac{r^{-1}}{1+|a_1|}+\dfrac{1}{1-r}\right)\left(\dfrac{|a_0|^2+|a_1|^2r^2}{|a_1|^2r^2}\right)-\dfrac{1}{|a_1|^2r^2},\; \mbox{for}\; k=1\vspace{2mm}\\
				\left(\dfrac{r^{-1}}{1+|a_1|}+\dfrac{1}{1-r}\right),\; \mbox{for}\; k\geq 2,
			\end{cases}
		\end{align*}
		we obtain (see \cite[Theorem 3(b)]{Ponnusamy-Vijayak-Wirths-RM-2020}):
		\begin{align*}
			\sum_{k=1}^{\infty}|a_k|r^k+\left(\frac{r^{-1}}{1+|a_1|}+\dfrac{1}{1-r}\right)\sum_{k=0}^{\infty}|a_k|^2r^{2k-1}\leq 1\; \mbox{for}\; r\leq \frac{5-\sqrt{17}}{2}.
		\end{align*}
		The number ${(5-\sqrt{17})}/{2}$ is sharp.\vspace{1.2mm}
		
		\item[(iv)] Suppose that $f\in\mathcal{B}$, $f(z)=\sum_{k=0}^{\infty}a_{qk}z^{qk}$ and $q\in\mathbb{N}$. Then for $\varphi_k(r)=r^{qk}\; (k\geq 0)$ and 
		\begin{align*}
			\Phi_k(r)=\left(\frac{1}{1+a}+\frac{r^q}{1-r^q}\right)\; (k\geq 1)
		\end{align*}
		Theorem \ref{B-thm-4.2} yields \cite[Corollary 1]{Liu-Liu-Ponnusamy-2021}, which actually improves upon \cite[Corollary 1]{Kayumov-Ponnusamy-JMAA-2018}, and also provides a generalized form of \cite[Theorem 2]{Ponnusamy-Vijayak-Wirths-RM-2020}.
		\begin{enumerate}
			\item[(a)] if $p=1$, then 
			\begin{align*}
				\sum_{k=0}^{\infty}|a_{qk}|r^{qk}+\left(\frac{1}{1+|a_0|}+\frac{r^q}{1-r^q}\right)\sum_{k=1}^{\infty}|a_{qk}|^2r^{2qk}\leq 1
			\end{align*}
			for $r\leq {1}/{\sqrt[q]{(2+a)}}$. The numbers ${1}/{\sqrt[q]{(2+a)}}$ and ${1}/{(1+a)}$ cannot be improved.
			\item[(b)] if $p=2$, then
			\begin{align*}
				a^2+\sum_{k=1}^{\infty}|a_{qk}|r^{qk}+\left(\frac{1}{1+a}+\frac{r^q}{1-r^q}\right)\sum_{k=1}^{\infty}|a_{qk}|^2r^{2qk}\leq 1
			\end{align*}
			for $r\leq {1}/{\sqrt[q]{2}}$. The numbers ${1}/{\sqrt[q]{2}}$ and ${1}/{(1+a)}$ cannot be improved.
		\end{enumerate}
		\vspace{3mm}
	\item[(II)] {\bf Improved Bohr inequalities}:  It is well-known that $S_r=S_r(f)$ denotes the planar integral
	\begin{align*}
		S_r=\int_{\mathbb D_r} |f'(z)|^2 d A(z).
	\end{align*}
	Note that if  $f(z)=\sum_{n=0}^\infty a_nz^n$, then $S_r=\pi \sum_{n=1}^\infty n|a_n|^2 r^{2n}.$  If $f$ is a univalent function, then $S_r$ is the area of  $f(\mathbb D_r)$.
	In the study of improved Bohr inequality, it is well-known that the role of sharp bounds of the coefficients $ a_n $ in the majorant series and sharp bounds 
	\begin{align*}
		\frac{S_r}{\pi}=\sum_{n=1}^{\infty}n|a_n|^2r^{2n}\leq \frac{r^2(1-|a_0|^2)^2}{(1-|a_0|^2r^2)^2} \,\,\,\,\mbox{for}\,\,\,\, 0<r\leq1/\sqrt{2}
	\end{align*}
	of the quantity $ S_r $ for the class of functions $f\in \mathcal{B} $ is established by Kayumov and Ponnusamy \cite{Kayumov-CRACAD-2018} are significant.\vspace{2mm}

		\item[(i)] For $\varphi_k(r)=r^k$ $(k\geq 0)$, we easily obtain (see \cite[Theorem 1]{Kayumov-CRACAD-2018}):
		\begin{enumerate}
			\item[(a)] if $p=1$ and $\Phi_k(r)=\frac{16}{9}k$, $(k\geq 1)$, then
			\begin{align*}
				\sum_{k=0}^{\infty}|a_k|r^k+\frac{16}{9}\left(\frac{S_r}{\pi}\right)\leq 1\; \mbox{for}\; r\leq \frac{1}{3}.
			\end{align*}
			The numbers $16/9$ and $1/3$ cannot be improved.\vspace{2mm}
			\item[(b)] if $p=2$ and $\Phi_k(r)=\frac{9}{8}k$, $(k\geq 1)$, then
			\begin{align*}
				a^2+\sum_{k=1}^{\infty}|a_k|r^k+\frac{9}{8}\left(\frac{S_r}{\pi}\right)\leq 1\; \mbox{for}\; r\leq \frac{1}{2}.
			\end{align*}
			The numbers $9/8$ and $1/2$ cannot be improved.\vspace{2mm}
		\end{enumerate}
		\item[(ii)] For $p=1$, $\varphi_k(r)=r^k\; (k\geq 0)$, and $\Phi_k(r)=\frac{1}{2r^k}\; (k\geq 1)$, $r\in (0, 1)$, we have (\cite[Theorem 2]{Kayumov-CRACAD-2018})
		\begin{align*}
			|a_0|+\sum_{k=1}^{\infty}\left(|a_k|+\frac{1}{2}|a_k|^2\right)r^k\leq 1\; \mbox{for}\; r\leq\frac{1}{3}.
		\end{align*}
		The numbers $1/2$ and $1/3$ cannot be improved.
		\item[(iii)]  For $\varphi_k(r)=r^k\; (k\geq 0)$,  we obtain (\cite[Theorem 4]{Liu-Liu-Ponnusamy-2021}) 
		\begin{enumerate}
			\item[(a)] if $p=1$ and 
			\begin{align*}
				\Phi_k(r)=\left(\frac{1}{1+a}+\frac{r}{1-r}\right)+\frac{8}{9}k,\; k\geq 1,
			\end{align*} 
			then we have 
			\begin{align*}
				\sum_{k=0}^{\infty}|a_k|r^k+\left(\frac{1}{1+a}+\frac{r}{1-r}\right)\sum_{k=1}^{\infty}|a_k|^2r^{2k}+\frac{8}{9}\left(\frac{S_r}{\pi}\right)\leq 1
			\end{align*}
			for $r\leq{1}/{3}$. 
			The numbers $1/3$ and $8/9$ cannot be improved. Moreover, \vspace{2mm}
			\item[(b)] if $p=2$ and 
			\begin{align*}
				\Phi_k(r)=\left(\frac{1}{1+a}+\frac{r}{1-r}\right)+\frac{9}{8}k,\; k\geq 1,
			\end{align*} 
			then we have 
			\begin{align*}
				a^2+\sum_{k=1}^{\infty}|a_k|r^k+\left(\frac{1}{1+a}+\frac{r}{1-r}\right)\sum_{k=1}^{\infty}|a_k|^2r^{2k}+\frac{9}{8}\left(\frac{S_r}{\pi}\right)\leq 1
			\end{align*}
			for $r\leq{1}/{(3-a)}$. 
			The numbers $1/(3-a)$ and $9/8$ cannot be improved.
		\end{enumerate}
\end{enumerate}
\subsection{Proof of Theorem \ref{B-thm-4.2}}
\begin{proof}[\bf Proof of Theorem \ref{B-thm-4.2}]
	For $f\in\mathcal{B}$, applying the Schwarz–Pick lemma yields the inequality $|a_k|\leq 1-|a|^2$ for $k\geq 1$, by a simple computation, we obtain
	\begin{align*}
		\mathcal{C}_f(\varphi, \Phi, p, r)&=|a_0|^p\varphi_0(r)+\sum_{k=1}^{\infty}|a_k|\varphi_k(r)+\sum_{k=1}^{\infty}\Phi_k(r)|a_k|^2\left(\varphi_k(r)\right)^2\\&\leq \varphi_0(r)-\left(\left(1-|a_0|^p\right)\right)\varphi_0(r)+\left(1-|a_0|^2\right)\sum_{k=1}^{\infty}\varphi_k(r)\\&\quad+\left(1-|a_0|^2\right)^2\sum_{k=1}^{\infty}\Phi_k(r)\left(\varphi_k(r)\right)^2\\&=\varphi_0(r)+\left(1-|a_0|^2\right)\bigg(\sum_{k=1}^{\infty}\varphi_k(r)-\left(\frac{1-|a_0|^p}{1-|a_0|^2}\right)\varphi_0(r)\\&\quad+\left(1-|a_0|^2\right)\sum_{k=1}^{\infty}\Phi_k(r)\left(\varphi_k(r)\right)^2\bigg).
	\end{align*}
	To further progress in the proof, we rely on the following inequality that has been proven in \cite{Kayu-Kham-Ponn-MJM-2021}
	\begin{align*}
		\frac{1-x^p}{1-x^2}\geq\frac{p}{2},\; \mbox{for all}\; x\in [0, 1]\; \mbox{and}\; p\in [0, p).
	\end{align*}
	Thus, we see that
	\begin{align*}
		\mathcal{C}_f(\varphi, \Phi, p, r)&\leq \varphi_0(r)+\left(1-|a_0|^2\right)\bigg(\sum_{k=1}^{\infty}\varphi_k(r)-\frac{p}{2}\varphi_0(r)+\left(1-|a_0|^2\right)\sum_{k=1}^{\infty}\Phi_k(r)\left(\varphi_k(r)\right)^2\bigg) \\&\leq \varphi_0(r)+\left(1-|a_0|^2\right)\Psi_p(r),
	\end{align*}
	where 
	\begin{align*}
		\Psi_p(r):=\sum_{k=1}^{\infty}\varphi_k(r)-\frac{p}{2}\varphi_0(r)+\left(1-|a_0|^2\right)\sum_{k=1}^{\infty}\Phi_k(r)\left(\varphi_k(r)\right)^2.
	\end{align*}
	Evidently, the desired inequality $ B_f(\varphi, p, r)\leq \varphi_0(r) $ holds if $ \Psi_p(r)\leq 0 $, i.e. if 
	\begin{align*}
		\varphi_0(r)\geq \frac{2}{p}\sum_{k=1}^{\infty}\left(\varphi_k(r)+\Phi_k(r)\left(\varphi_k(r)\right)^2\right).
	\end{align*}
	Thus the inequality of the theorem is established for $ |z|=r\leq R $, where $ R $ is the minimal positive root of equation
	\begin{align*}
		\varphi_0(r)=\frac{2}{p}\sum_{k=1}^{\infty}\left(\varphi_k(r)+\Phi_k(r)\left(\varphi_k(r)\right)^2\right).
	\end{align*} \vspace{1.2mm}
	
	To complete the proof, it is sufficient to show that the result is sharp. In order to show that, we consider the function $ f_a $ which is defined as
	\begin{align*}
		f_a(z)=\frac{z-a}{1-az}=a+\sum_{k=1}^{\infty}\left(1-a^2\right)a^{k-1}z^k.
	\end{align*}
	A simple computation gives that
	\begin{align*}
		&|a_0|^p\varphi_0(r)+\sum_{k=1}^{\infty}|a_k|\varphi_k(r)+\Phi_k(r)\sum_{k=1}^{\infty}|a_k|^2\left(\varphi_k(r)\right)^2\\&=\varphi_0(r)+\left(1-a\right)\bigg((1+a)\sum_{k=1}^{\infty}a^{k-1}\varphi_k(r)+\Phi_k(r)(1+a)\left(1-a^2\right)\sum_{k=1}^{\infty}a^{2k-2}\left(\varphi_k(r)\right)^2\\&\quad-\frac{1-a^p}{1-a}\varphi_0(r)\bigg)\\&=\varphi_0(r)+\left(1-a\right)\bigg(2\sum_{k=1}^{\infty}a^{k-1}\varphi_k(r)+\sum_{k=1}^{\infty}a^{2k-2}\left(\varphi_k(r)\right)^2-p\varphi_0(r)\bigg)\\&\quad+\left(1-a^2\right)^2\sum_{k=1}^{\infty}a^{k-1}\varphi_k(r)+\left(\Phi_k(r)\left(1-a^2\right)^2-2(1-a)\right)\sum_{k=1}^{\infty}a^{2k-2}\left(\varphi_k(r)\right)^2\\&\quad-\left((1-a^p)-p(1-a)\right)\varphi_0(r).
	\end{align*}
	Then for $ a\rightarrow 1^{-} $, it is easy to see that 
	\begin{align*}
		&|a_0|^p\varphi_0(r)+\sum_{k=1}^{\infty}|a_k|\varphi_k(r)+\Phi_k(r)\sum_{k=1}^{\infty}|a_k|^2\left(\varphi_k(r)\right)^2\\&=\varphi_0(r)+(1-a)\left(2\sum_{k=1}^{\infty}\bigg\{\varphi_k(r)+\Phi_k(r)\left(\varphi_k(r)\right)^2\bigg\}-p\varphi_0(r)\right)+O\left(\left(1-a\right)^2\right)\\&>\varphi_0(r)
	\end{align*}
	if 
	\begin{align*}
		2\sum_{k=1}^{\infty}\bigg\{\varphi_k(r)+\Phi_k(r)\left(\varphi_k(r)\right)^2\bigg\}-p\varphi_0(r)>0
	\end{align*}
	\textit{i.e.,} if 
	\begin{align*}
		\varphi_0(r)<\frac{2}{p}\sum_{k=1}^{\infty}\bigg\{\varphi_k(r)+\Phi_k(r)\left(\varphi_k(r)\right)^2\bigg\}.
	\end{align*}
	This establishes the sharpness and the proof is completed.
\end{proof}
\section{Generalized Bohr inequalities for certain class $\tilde{G}^0_{\mathcal{H}}(\beta)$ of harmonic mappings}
Methods of harmonic mappings have been applied to study and solve the fluid flow problems (see \cite{Aleman-2012,Constantin-2017}). For example, in 2012, Aleman and Constantin \cite{Aleman-2012} established a connection between harmonic mappings and ideal fluid flows. In fact, Aleman and Constantin have developed an ingenious technique to solve the incompressible two dimensional Euler equations in terms of univalent harmonic mappings (see \cite{Constantin-2017} for details).\vspace{1.2mm}

A complex-valued function $ f(z)=u(x,y)+iv(x,y) $ is called harmonic in $ U $ if both $ u $ and $ v $ satisfy the Laplace's equation $ \bigtriangledown^2 u=0 $ and $ \bigtriangledown^2 v=0 $, where
\begin{equation*}
	\bigtriangledown^2:=\frac{\partial^2}{\partial x^2}+\frac{\partial^2}{\partial y^2}.
\end{equation*}
\par It is well-known that under the assumption $ g(0)=0 $, the harmonic function $ f $ has the unique canonical representation $ f=h+\overline{g} $, where $ h $ and $ g $ are analytic functions in $ U $, respectively called, analytic and co-analytic parts of $ f $.  If in addition $ f $ is univalent then we say that $ f $ is univalent harmonic on a domain $ \Omega $. A locally univalent harmonic mapping $ f=h+\overline{g} $ is sense-preserving whenever its Jacobian $ J_f(z):=|f_{z}(z)|^2-|f_{\bar{z}}(z)|=|h^{\prime}(z)|^2-|g^{\prime}(z)|^2>0 $ for $ z\in \mathbb{D} $.\vspace{1.2mm} 

We note that the inequality $ M_f(r)\leq 1 $ can be written in distance formulations as follows
\begin{align*}
	\sum_{n=1}^{\infty}|a_n|r^n\leq 1-|a_0|=1-|f(0)|=d(f(0), \partial\mathbb{D}),
\end{align*}
where $ d $ is the Euclidean distance between $ f(0) $ and the boundary $ \partial\mathbb{D} $ of the unit disk $ \mathbb{D} $. In view of this distance formulations, Bohr's phenomenon has been studied recently (see \cite{Aha-Allu-RMJ-2022,Aha-CMFT-2022,Arora-CVEE-2022}). 
Inspired by the notion of Rogosinski’s inequality and Rogosinski’s radius investigated in \cite{Rogosinski-1923}, Kayumov \emph{et al} \cite{Kayumov-Khammatova-Ponnusamy-JMAA-2021} obtained the following Bohr-Rogosinski inequality and Bohr-Rogosinski radius for the class $ \mathcal{B}.$
\begin{thm}\cite{Kayumov-Khammatova-Ponnusamy-JMAA-2021}\label{B-th-33.11}
	Suppose that $ f(z)=\sum_{n=0}^{\infty}a_nz^n\in\mathcal{B} $. Then 
	\begin{equation*}
		|f(z)|+\sum_{n=N}^{\infty}|a_n|r^n\leq 1\;\;\mbox{for}\;\; r\leq R_N,
	\end{equation*}
	where $ R_N $ is the unique root of the equation $ 2(1+r)r^N-(1-r)^2=0 $ in the interval $(0,1)$. The radius $ R_N $ is the best possible. Moreover, 
	\begin{equation*}
		|f(z)|^2+\sum_{n=N}^{\infty}|a_n|r^n\leq 1\;\;\mbox{for}\;\; r\leq R^{\prime}_N,
	\end{equation*}
	where $ R^{\prime}_N $ is the unique root of the equation $ (1+r)r^N-(1-r)^2=0 $ in the interval $(0,1)$. The radius $ R^{\prime}_N $ is the best possible.
\end{thm}
In what follows, $ \lfloor x \rfloor $ denotes the largest integer not more than $ x $, where $ x $ is a real number. It is pertinent to mention a key result from \cite{Liu-Liu-Ponnusamy-2021} that has been utilized to establish a more refined form of the Bohr-type inequalities.
\begin{lem}\cite{Liu-Liu-Ponnusamy-2021}\label{Lemm-1.2}
	Suppose that $ f(z)=\sum_{n=0}^{\infty}a_nz^n\in \mathcal{B} $. Then for any $N\in\mathbb{N}$, the following inequality holds:
	\begin{align*}
		\sum_{n=N}^{\infty}|a_n|r^n&+sgn(t)\sum_{n=1}^{t}|a_n|^2\dfrac{r^N}{1-r}+\left(\dfrac{1}{1+|a_0|}+\dfrac{r}{1-r}\right)\sum_{n=t+1}^{\infty}|a_n|^2r^{2n}\\&\leq (1-|a_0|^2)\dfrac{r^N}{1-r}, 
	\end{align*}
	\;\mbox{for}\; $|z|= r\in[0,1),$ where $t=\lfloor{(N-1)/2}\rfloor$.
\end{lem}
By applying Lemma \ref{Lemm-1.2}, the next result is obtained as refined versions of classical Bohr inequality Liu \emph{et al.} in \cite{Liu-Liu-Ponnusamy-2021}.
\begin{thm}\cite{Liu-Liu-Ponnusamy-2021}\label{Th-4.2}
	Suppose that $ f\in\mathcal{B} $ and $ f(z)=\sum_{n=0}^{\infty}a_nz^n $. For $ n\in\mathbb{N} $, let $ t=\lfloor (N-1)/2 \rfloor $. Then 
	\begin{align*}
		|f(z)|+\sum_{n=N}^{\infty}|a_n|r^n+sgn(t)\sum_{n=1}^{t}|a_n|^2\frac{r^N}{1-r}+\left(\frac{1}{1+|a_0|}+\frac{r}{1-r}\right)\sum_{n=t+1}^{\infty}|a_n|^2r^{2n}\leq 1
	\end{align*}
	for $ |z|=r\leq R_N $, where $ R_N $ is as in Theorem \ref{B-th-33.11}. The radius $ R_N $ is best possible. Moreover, 
	\begin{align*}
		|f(z)|^2+\sum_{n=N}^{\infty}|a_n|r^n+sgn(t)\sum_{n=1}^{t}|a_n|^2\frac{r^N}{1-r}+\left(\frac{1}{1+|a_0|}+\frac{r}{1-r}\right)\sum_{n=t+1}^{\infty}|a_n|^2r^{2n}\leq 1
	\end{align*}
	for $ |z|=r\leq R^{\prime}_N $, where $ R^{\prime}_N $ is as in Theorem \ref{B-th-33.11}. The radius $ R_N $ is best possible.
\end{thm}
For the case $ N=1 $, it is clear that $ R_1=\sqrt{5}-2 $ and $ R^{\prime}_1=1/3 $. Moreover, the following result is established showing that the two constants $ R_1=\sqrt{5}-2 $ and $ R^{\prime}_1=1/3 $ can be improved for individual functions in $ \mathcal{B} $ in the context of Theorem \ref{B-th-33.11} and Theorem \ref{Th-4.2} with $ N=1 $. 
\begin{thm}\cite{Liu-Liu-Ponnusamy-2021}\label{Th-4.3}
	Suppose that $ f\in\mathcal{B} $ and $ f(z)=\sum_{n=0}^{\infty}a_nz^n $. Then 
	\begin{align*}
		A(z):=|f(z)|+\sum_{n=1}^{\infty}|a_n|r^n+\left(\frac{1}{1+|a_0|}+\frac{r}{1-r}\right)\sum_{n=1}^{\infty}|a_n|^2r^{2n}\leq 1
	\end{align*}
	for $ |z|=r\leq r_{a_0}=2/(3+|a_0|+\sqrt{5}(1+|a_0|)) $. The radius $ r_{a_0}$ is best possible and $ r_{a_0}>\sqrt{5}-2 $. Moreover, 
	\begin{align*}
		B(z):=|f(z)|^2+\sum_{n=1}^{\infty}|a_n|r^n+\left(\frac{1}{1+|a_0|}+\frac{r}{1-r}\right)\sum_{n=1}^{\infty}|a_n|^2r^{2n}\leq 1
	\end{align*}
	for $ |z|=r\leq r^{\prime}_{a_0} $, where $ r^{\prime}_{a_0} $ is unique positive root of the equation 
	\begin{align*}
		\left(1-|a_0|^3\right)r^3-(1+2|a_0|)r^2-2r+1=0.
	\end{align*} The radius $ r^{\prime}_{a_0} $ is best possible. Further, $ 1/3<r^{\prime}_{a_0}<1/(2+|a_0|) $. 
\end{thm}
The Bohr phenomenon can be generalized (see \cite{Allu-Hal-BSM-2021}) as: for a given domain $ G\subseteq\mathbb{C} $, to find the largest radius $ r_G>0 $ such that 
\begin{align*}
	\sum_{n=1}^{\infty}|a_n|r^n\leq d(f(0), \partial G)
\end{align*}
for all functions $ f $ belong to the class of analytic functions in $ \mathbb{D} $ such that $ f(\mathbb{D}\subseteq G) $. \vspace{1.2mm}
Let $ \mathcal{H} $ be the class of all complex-valued harmonic mappings of the form Mf=h+$\overline{g}$ defined on unit disk $ \mathbb{D} $, where $ h $ and $ g $ are analytic functions on $ \mathbb{D} $ with the normalization $ h(0)=0=h^{\prime}(0)-1 $ and $ g(0)=0 $. Set 
\begin{align*}
	\mathcal{H}_0=\{f=h+\overline{g}\in\mathcal{H} : g^{\partial}(0)=0\}.
\end{align*}
According to the study carried out by Lewy \cite{Lewy-BAMS-1936}, a function $ f=h+\overline{g}\in\mathcal{H} $ is locally univalent
and sense-preserving on $ \mathbb{D} $ if and only if its Jacobian $ J_f(z) $ is positive in $ \mathbb{D} $, where 
\begin{align*}
	J_f(z)=|f_{z}(z)|^2-|f_{\overline{z}}(z)|^2=|h^{\prime}(z)|^2-|g^{\prime}(z)|^2.
\end{align*}
In light of this finding, it is evident that $ J_f(z)>0 $ in $ \mathbb{D} $ if, and only if, $ h^{\prime}(z)\neq 0 $ in $ \mathbb{D} $ and the (second complex) dilation $ \omega(z)=g^{\prime}(z)/h^{\prime}(z) $ of $ f=h+\overline{g} $ is analytic in $ \mathbb{D} $ and has the property that $ |\omega(z)|<1 $ for $ z\in\mathbb{D} $.\vspace{2mm}

In $ 2016 $, Liu and Ponnusamy \cite{Li-Ponn-CMJ-2016} have studied the class $ \tilde{G}^0_{\mathcal{H}}(\beta) $, where
\begin{align*}
	\tilde{G}^0_{\mathcal{H}}(\beta)=\bigg\{f=h+\bar{g}\in\mathcal{H}_0 : \rm Re\left(\frac{h(z)}{z}\right)-\beta>\bigg|\frac{g(z)}{z}\bigg|\; \mbox{for}\; z\in\mathbb{D}\bigg\}
\end{align*}
and established that if $ f\in \tilde{G}^0_{\mathcal{H}}(\beta) $, then the harmonic convolution is a univalent and close-to-convex harmonic
function in the unit disk provided certain conditions for parameter $ \beta $ is satisfied. It can be shown that $ \tilde{G}^0_{\mathcal{H}}(\beta)\subseteq \tilde{G}^0_{\mathcal{H}}(0) $ for $ 0\leq \beta<1 $.\vspace{1.2mm}

 The Bohr's phenomenon for certain class $ \tilde{G}^0_{\mathcal{H}}(\beta) $ consisting of functions $ f $ of the form 
\begin{align}\label{B-eq-3.1}
	f(z)=h(z)+\overline{g(z)}=z+\sum_{n=2}^{\infty}a_nz^n+\overline{\sum_{n=2}^{\infty}b_nz^n}
\end{align}
is as follows: find the largest radius $ r_f\in (0, 1) $ such that the following inequality 
\begin{align}\label{B-eqq-3.2}
	|z|+\sum_{n=2}^{\infty}(|a_n|+|b_n|)|z|^n\leq d(f(0), \partial f(\mathbb{D}))
\end{align}
holds for all $ |z|=r\leq r_f $, and for all $ f\in\tilde{G}^0_{\mathcal{H}}(\beta) $. The radius $ r_f $ is called Bohr's radius for the class $ \tilde{G}^0_{\mathcal{H}}(\beta) $. The sharp coefficient bound $ |a_n|+|b_n|\leq 2(1-\beta) $ for the class $ \tilde{G}^0_{\mathcal{H}}(\beta) $ is obtained in \cite{Allu-Hal-BSM-2021}
and using this the following result is established.
\begin{thm}\cite{Allu-Hal-BSM-2021}\label{B-theorem-3.1}
	Let $ f=h+\bar{g}\in \tilde{G}^0_{\mathcal{H}}(\beta) $ for $ 0<\beta<1/2 $ be given by \eqref{B-eq-3.1}. Then the inequality \eqref{B-eqq-3.2} is satisfied for $ |z|=r\leq r_f:=\left(-1-\beta+\sqrt{1+6\beta-7\beta^2}\right)/2(1-2\beta) $, where $ r_f $ is the unique root in $ (0, 1) $ of the equation $ (1-2\beta)r^2+(1+\beta)r-\beta=0 $. The constant $ r_f $ is best possible.
\end{thm}

In $ 2010 $, Abu-Muhanna \cite{Abu-CVEE-2010} considered first time the Bohr radius for the class of complex-valued harmonic function $ f=h+\bar{g} $ defined in $\mathbb{D}$ with $|f(z)|<1$ for all $z\in \mathbb{D}.$ After this, Kayumov \emph{et al.}\cite{Kayumov-Ponnusamy-Shakirov-MN-2017} studied Bohr radius for locally univalent harmonic mappings, Liu and Ponnusamy \cite{Liu-Ponnusamy-BMMSS-2019} have determined the Bohr radius for $k$-quasiconformal  harmonic mappings, Evdoridis \emph{et al.}\cite{Evd-Pon-Ras-Antti-IM-2019} studied an improved version of  the Bohr's inequality for locally univalent harmonic mappings, Ahamed \cite{Aha-CMFT-2022,MBA-CVEE-2022}, Ahamed and Allu \cite{Aha-Allu-CMB-2023} have studied refined Bohr-Rogosinski inequalities for certain classes of harmonic mappings. Recently, Arora \cite{Arora-CVEE-2022} have studied Bohr-type inequality for harmonic functions with lacunary series in complex Banach space. \vspace{1.2mm} 

In light of the findings in \cite{Liu-Liu-Ponnusamy-2021}, our objective in this paper is twofold: firstly, to obtain an improved version of Theorem \ref{B-theorem-3.1}, and secondly, to establish a harmonic analogue of the Bohr inequalities for bounded analytic functions on $ \mathbb{D} $. Henceforth, for $ f=h+\bar{g}\in \tilde{G}^0_{\mathcal{H}}(\beta) $ and $ |z|=r $, we first denote the functional
\begin{align*}
	S^f_{\mu, \lambda, m, N}(r):&=|f(z)|^m+\sum_{n=N}^{\infty}\left(|a_n|+|b_n|\right)r^n+\mu sgn(t)\sum_{n=1}^{t}\left(|a_n|+|b_n|\right)^2\frac{r^N}{1-r}\\&\quad+\lambda\left(1+\frac{r}{1-r}\right)\sum_{n=t+1}^{\infty}\left(|a_n|+|b_n|\right)^2r^{2n}.
\end{align*}
\begin{thm}\label{B-thm-3.1}
	Let $ f=h+\bar{g}\in \tilde{G}^0_{\mathcal{H}}(\beta) $ for $ 0<\beta<1/2 $ be given by \eqref{B-eq-3.1}. Then the inequality $ S^f_{\mu, \lambda, m, N}(r)\leq d(f(0), \partial f(\mathbb{D})) $ is satisfied for for $ |z|=r\leq R^{m, N, t}_{\mu, \lambda, \beta}(r) $, where $ R^{m, N, t}_{\mu, \lambda, \beta}(r) $ is the unique root in $ (0, 1) $ of equation
	\begin{align*}
		&\left(\beta r+(1-\beta)\left(\dfrac{1+r}{1-r}\right)r\right)^m+\dfrac{2(1-\beta)r^N}{1-r}+4\mu sgn(t)\frac{ (1-\beta)^2r^Nt}{1-r}\\&\quad+\left(1+\frac{r}{1-r}\right)\frac{4\lambda (1-\beta)^2r^{2t+2}}{1-r^2}-\beta=0.
	\end{align*} The constant $ R^{m, N, t}_{\mu, \lambda, \beta}(r) $ is best possible.
\end{thm}
\subsection{Corollaries of Theorem \ref{B-thm-3.1}}
Due to the implications of Theorem \ref{B-theorem-3.1}, we derive the subsequent corollary, which provides a harmonic counterpart to Theorem \ref{B-th-33.11} specifically for the class $ \tilde{G}^0_{\mathcal{H}}(\beta) $, particularly when $ m=1, 2 $.
\begin{cor}
	Let $ f=h+\bar{g}\in \tilde{G}^0_{\mathcal{H}}(\beta) $ for $ 0<\beta<1/2 $ be given by \eqref{B-eq-3.1}. Then the inequality 
	\begin{align*}
		|f(z)|^m+\sum_{n=N}^{\infty}\left(|a_n|+|b_n|\right)r^n+	\leq d(f(0), \partial f(\mathbb{D}))
	\end{align*} is satisfied for for $ |z|=r\leq R^{m, N}_{\beta}(r) $, where $ R^{m, N}_{\beta}(r) $ is the unique root in $ (0, 1) $ of equation
	\begin{align*}
		\left(\beta r+(1-\beta)\left(\dfrac{1+r}{1-r}\right)r\right)^m+\dfrac{2(1-\beta)r^N}{1-r}-\beta=0.
	\end{align*} The constant $ R^{m, N}_{\beta}(r) $ is best possible.
\end{cor}
In the special case where $ m=1, 2 $ and $ \mu=1=\lambda $, the subsequent result establishes the Bohr inequality represents a harmonic counterpart to the inequalities presented in Theorem \ref{Th-4.2} (see \cite[Theorem 1, p. 7]{Liu-Liu-Ponnusamy-2021}).
\begin{cor} Let $ f=h+\bar{g}\in \tilde{G}^0_{\mathcal{H}}(\beta) $ for $ 0<\beta<1/2 $ be given by \eqref{B-eq-3.1}. Then
	\begin{align*}
		&|f(z)|+\sum_{n=N}^{\infty}\left(|a_n|+|b_n|\right)r^n+ sgn(t)\sum_{n=1}^{t}\left(|a_n|+|b_n|\right)^2\frac{r^N}{1-r}\\&\quad+\left(1+\frac{r}{1-r}\right)\sum_{n=t+1}^{\infty}\left(|a_n|+|b_n|\right)^2r^{2n}d(f(0), \partial f(\mathbb{D}))
	\end{align*}
	is satisfied for for $ |z|=r\leq R^{1, N, t}_{1, 1, \beta}(r) $, where $ R^{1, N, t}_{1, 1, \beta}(r) $ is the unique root in $ (0, 1) $ of equation 
	\begin{align*}
		&\beta r+(1-\beta)\left(\dfrac{1+r}{1-r}\right)r+\dfrac{2(1-\beta)r^N}{1-r}+4\; sgn(t)\frac{ (1-\beta)^2r^Nt}{1-r}\\&\quad+\left(1+\frac{r}{1-r}\right)\frac{4 (1-\beta)^2r^{2t+2}}{1-r^2}-\beta=0.
	\end{align*} 
	The constant $ R^{1, N, t}_{1, 1, \beta}(r) $ is best possible. Moreover,
	\begin{align*}
		&|f(z)|^2+\sum_{n=N}^{\infty}\left(|a_n|+|b_n|\right)r^n+ sgn(t)\sum_{n=1}^{t}\left(|a_n|+|b_n|\right)^2\frac{r^N}{1-r}\\&\quad+\left(1+\frac{r}{1-r}\right)\sum_{n=t+1}^{\infty}\left(|a_n|+|b_n|\right)^2r^{2n}\leq d(f(0), \partial f(\mathbb{D}))
	\end{align*}
	is satisfied for for $ |z|=r\leq R^{2, N, t}_{1, 1, \beta}(r) $, where $ R^{2, N, t}_{1, 1, \beta}(r) $ is the unique root in $ (0, 1) $ of equation	
	\begin{align*}
		&\left(\beta r+(1-\beta)\left(\dfrac{1+r}{1-r}\right)r\right)^2+\dfrac{2(1-\beta)r^N}{1-r}+4 sgn(t)\frac{ (1-\beta)^2r^Nt}{1-r}\\&\quad+\left(1+\frac{r}{1-r}\right)\frac{4(1-\beta)^2r^{2t+2}}{1-r^2}-\beta=0.
	\end{align*} The constant $ R^{2, N, t}_{1, 1, \beta}(r) $ is best possible.
\end{cor}
When $ m=1, 2 $ and $ \mu $ is allowed to vary, with $ \lambda=1 $, the following result establishes that the Bohr inequality can be seen as a harmonic analog to the inequalities described in Theorem \ref{Th-4.3} (as shown in \cite[Theorem 2, p. 7]{Liu-Liu-Ponnusamy-2021}).
\begin{cor} Let $ f=h+\bar{g}\in \tilde{G}^0_{\mathcal{H}}(\beta) $ for $ 0<\beta<1/2 $ be given by \eqref{B-eq-3.1}. Then
	\begin{align*}
		&|f(z)|+\sum_{n=1}^{\infty}\left(|a_n|+|b_n|\right)r^n+\left(1+\frac{r}{1-r}\right)\sum_{n=1}^{\infty}\left(|a_n|+|b_n|\right)^2r^{2n}\leq d(f(0), \partial f(\mathbb{D}))
	\end{align*}
	is satisfied for for $ |z|=r\leq R^{1, 1, 0}_{0, 1, \beta}(r) $, where $ R^{1, 1, 0}_{0, 1, \beta}(r) $ is the unique root in $ (0, 1) $ of equation
	\begin{align*}
		&\beta r+(1-\beta)\left(\dfrac{1+r}{1-r}\right)r+\dfrac{2(1-\beta)r^N}{1-r}\\&\quad+\left(1+\frac{r}{1-r}\right)\frac{4 (1-\beta)^2r^{2}}{1-r^2}-\beta=0.
	\end{align*} The constant $ R^{1, 1, 0}_{0, 1, \beta}(r) $ is best possible. Moreover,
	\begin{align*}
		&|f(z)|^2+\sum_{n=1}^{\infty}\left(|a_n|+|b_n|\right)r^n+\left(1+\frac{r}{1-r}\right)\sum_{n=1}^{\infty}\left(|a_n|+|b_n|\right)^2r^{2n}\leq d(f(0), \partial f(\mathbb{D}))
	\end{align*}
	is satisfied for for $ |z|=r\leq R^{2, 1, 0}_{0, 1, \beta}(r) $, where $ R^{2, 1, 0}_{0, 1, \beta}(r)  $ is the unique root in $ (0, 1) $ of equation
	\begin{align*}
		&\left(\beta r+(1-\beta)\left(\dfrac{1+r}{1-r}\right)r\right)^2+\dfrac{2(1-\beta)r^N}{1-r}\\&\quad+\left(1+\frac{r}{1-r}\right)\frac{4 (1-\beta)^2r^{2}}{1-r^2}-\beta=0.
	\end{align*} 
	The constant $ R^{2, 1, 0}_{0, 1, \beta}(r)  $ is best possible.
\end{cor}
For further improvement of Bohr's inequality for harmonic mappings for the class $ \tilde{G}^0_{\mathcal{H}}(\beta) $, it is natural to raise the following question.
\begin{ques}\label{Q-3.1}
	In view of the general setting proposed in \cite{Kayu-Kham-Ponn-MJM-2021} with a change of basis from $\{r^n\}_{n=1}^{\infty}$ to $\{\varphi_n(r)\}_{n=1}^{\infty}$, is it possible to present a general result improving Theorem \ref{B-theorem-3.1}?
\end{ques}
As a response to Question \ref{Q-3.1}, we obtained the following result.
\begin{thm}\label{Th-3.6}
Let $ f=h+\bar{g}\in \tilde{G}^0_{\mathcal{H}}(\beta) $ for $ 0\leq\beta<1 $ be given by \eqref{B-eq-3.1}. If $ \{\varphi_n(r)\}_{n=1}^{\infty} $ is an increasing sequence of non-negative continuous functions in $ [0, 1) $ such that the series 
\begin{align*}
	\varphi_0(r)+\sum_{n=1}^{\infty}\varphi_n(r)
\end{align*} 
converges locally uniformly with respect to $ r\in [0, 1) $ and 
\begin{align}\label{Eq-3.3}
	\sum_{n=2}^{\infty}\varphi_n(0)<\frac{\beta}{2(1-\beta)},
\end{align} then the inequality 
\begin{align}
	r\varphi_1(r)+\sum_{n=2}^{\infty}(|a_n|+|b_n|)\varphi_n(r)\leq d(f(0), \partial f(\mathbb{D}))
\end{align} is satisfied for $ |z|=r\leq R_{\beta}$, where $ R_{\beta} $ is the unique root in $ (0, 1) $ of the equation $F_{\beta}(r)=0$, where
\begin{align}\label{Eq-3.4}
F_{\beta}(r)=r\varphi_0(r)+2(1-\beta)\sum_{n=2}^{\infty}\varphi_n(r)-\beta.
\end{align} 
The constant $ R_{\beta} $ is best possible.
\end{thm}
\begin{rem}
	The conclusion reached in Theorem \ref{Th-3.6} exposes the following:
	\begin{enumerate}
		\item[(i)] In the case where $\varphi_1(r)=1$ and $\varphi_n(r)=r^n$ $(n\geq 2)$ are chosen with $0<\beta<1/2$, it is easy to see that $\sum_{n=2}^{\infty}\varphi_n(0)=0<{\beta}/{2(1-\beta)}$ is satisfied trivially, hence Theorem \ref{Th-3.6} precisely reduces to Theorem \ref{B-theorem-3.1}. Consequently, Theorem \ref{Th-3.6} provides a broader formulation of Theorem \ref{B-theorem-3.1}, thereby successfully addressing Question \ref{Q-3.1}.
		\item[(ii)] By considering Theorem \ref{Th-3.6}, which extends the range of $\beta$ from $(0, 1/2)$ to $[0, 1)$, the applicability of Theorem \ref{B-theorem-3.1} has been broadened.
	\end{enumerate}
\end{rem}
\subsection{Some applications of Theorem \ref{Th-3.6}}
In the range $ 0\leq\beta<1 $ and with suitable selections of the sequence $ {\varphi_n(r)} $ in Theorem \ref{Th-3.6}, as a consequence, we arrived at the subsequent immediate results.
\begin{cor}
	Let $ f=h+\bar{g}\in \tilde{G}^0_{\mathcal{H}}(\beta) $ for $ 0\leq\beta<1 $ be given by \eqref{B-eq-3.1}. Let $ \varphi_1(r)=1 $ and $ \varphi_n(r)=(n+1)r^n $ for $ n\geq 2 $. Then, we have 
	\begin{align*}
		r+\sum_{n=2}^{\infty}(|a_n|+|b_n|)(n+1)r^n\leq d(f(0), \partial f(\mathbb{D}))
	\end{align*}
	for $ |z|=r\leq R_1(\beta) $, where $ R_1(\beta) $ is the unique root in $ (0, 1) $ of the equation 
	\begin{align*}
	r+2(1-\beta)\left(\frac{(3-2r)r^2}{(1-r)^2}\right)-\beta=0.
	\end{align*}
	The radius $ R_1(\beta) $ is best possible.
\end{cor}
By a straightforward calculation, it can be easily shown that for $ 0<r<1 $, 
\begin{align*}
	\sum_{n=2}^{\infty}nr^n=\frac{(2-r)r}{(1-r)^2},\; \sum_{n=2}^{\infty}n^2r^n=\frac{(4-3r+r^2)r^2}{(1-r)^3},\; \sum_{n=2}^{\infty}n^3r^n=\frac{(8-5r+4r^2-r^3)r^2}{(1-r)^4}.
\end{align*}
Taking into consideration these estimates, Theorem \ref{Th-3.6} leads to the following set of conclusions.
\begin{cor}
	Let $ f=h+\bar{g}\in \tilde{G}^0_{\mathcal{H}}(\beta) $ for $ 0\leq\beta<1 $ be given by \eqref{B-eq-3.1}. Let $ \varphi_1(r)=1 $ and $ \varphi_n(r)=n^{\alpha}r^n $ for $ n\geq 2 $. Then
	\begin{enumerate}
		\item[(i)] for $ \alpha=1 $, we have $	r+\sum_{n=2}^{\infty}(|a_n|+|b_n|)nr^n\leq d(f(0), \partial f(\mathbb{D}))$ for $ |z|=r\leq R_2(\beta) $, where $ R_2(\beta) $ is the unique root in $ (0, 1) $ of the equation 
		\begin{align*}
			r+2(1-\beta)\left(\frac{(2-r)r}{(1-r)^2}\right)-\beta=0.
		\end{align*}
		The radius $ R_2(\beta) $ is best possible. Moreover, \vspace{1.2mm}
		\item[(ii)] for $ \alpha=2 $, we have $	r+\sum_{n=2}^{\infty}(|a_n|+|b_n|)n^2r^n\leq d(f(0), \partial f(\mathbb{D}))$ for $ |z|=r\leq R_3(\beta) $, where $ R_3(\beta) $ is the unique root in $ (0, 1) $ of the equation 
		\begin{align*}
			r+2(1-\beta)\left(\frac{(4-3r+r^2)r^2}{(1-r)^3}\right)-\beta=0.
		\end{align*}
		The radius $ R_3(\beta) $ is best possible.
		\item[(iii)] for $ \alpha=3 $, we have $r+\sum_{n=2}^{\infty}(|a_n|+|b_n|)n^3r^n\leq d(f(0), \partial f(\mathbb{D}))$ for $ |z|=r\leq R_4(\beta) $, where $ R_4(\beta) $ is the unique root in $ (0, 1) $ of the equation 
		\begin{align*}
			r+2(1-\beta)\left(\frac{(8-5r+4r^2-r^3)r^2}{(1-r)^4}\right)-\beta=0.
		\end{align*}
		The radius $ R_4(\beta) $ is best possible.
	\end{enumerate}
\end{cor}
\begin{table}[ht]
	\centering
	\begin{tabular}{|l|l|l|l|l|l|l|l|l|l|}
		\hline
		$\;\;\beta$& $\;\;0.1$&$\;\;0.2$& $\;\;0.3$& $\;\; 0.4 $& $\;\;0.5$&$\;\;0.6$& $\;\;0.7$& $ \;\;0.8 $& $ \;\;0.9 $ \\
		\hline
		$R_1(\beta)$& $0.070 $&$0.119 $& $0.161$& $0.202$& $0.242 $&$0.287 $& $0.338$& $0.402$ &$ 0.5 $\\
		\hline
		$R_2(\beta)$& $0.021 $&$0.045 $& $0.072$& $0.104$& $0.142 $&$0.188 $& $0.245$& $0.324$ &$ 0.447 $\\
		\hline
		$R_3(\beta)$& $0.064 $&$0.106 $& $0.141$& $0.174$& $0.207 $&$0.242 $& $0.282$& $0.333$ &$ 0.412 $\\
		\hline
		$R_4(\beta)$& $0.052 $&$0.083 $& $0.108$& $0.131$& $0.155 $&$0.181 $& $0.211$& $0.249$ &$ 0.312 $\\
		\hline
	\end{tabular}
	\vspace{2.5mm}
	\caption{This table exhibits the the approximate values of the roots $ R_1(\beta) $, $ R_2(\beta) $, $ R_3(\beta) $, and $ R_4(\beta) $ for different values of $ \beta\in [0, 1) $.}
	\label{tabel-1}
\end{table}
\subsection{Proof of Theorems \ref{B-thm-3.1} and \ref{Th-3.6}}
\begin{proof}[\bf Proof of Theorem \ref{B-thm-3.1}]
	Let $ f=h+\bar{g}\in \tilde{G}^0_{\mathcal{H}}(\beta) $ for $ 0\leq\beta<1 $ be given by \eqref{B-eq-3.1}. A simple computation shows that
	\begin{align}
		\begin{cases}
			|f(z)|\leq \beta|z|+(1-\beta)\left(\dfrac{1+|z|}{1-|z|}\right)|z|=\beta r+(1-\beta)\left(\dfrac{1+r}{1-r}\right)r,\vspace{2mm}\\
			\displaystyle\sum_{n=N}^{\infty}\left(|a_n|+|b_n|\right)r^n=2(1-\beta)\sum_{n=N}^{\infty}r^n=\dfrac{2(1-\beta)r^N}{1-r},\vspace{2mm}\\
			\displaystyle\sum_{n=t+1}^{\infty}\left(|a_n|+|b_n|\right)^2r^{2n}\leq 4(1-\beta)^2\sum_{n=t+1}^{\infty}r^{2n}=\frac{4(1-\beta)^2r^{2t+2}}{1-r^2},\vspace{2mm}\\
			\displaystyle \sum_{n=1}^{t}\left(|a_n|+|b_n|\right)^2\dfrac{r^N}{1-r}\leq \frac{4(1-\beta)^2r^Nt}{1-r}.
		\end{cases}
	\end{align}
	Thus, it is easy to see that
	\begin{align*}
		S^f_{\mu, \lambda, m, N}(r)\leq & \left(\beta r+(1-\beta)\left(\dfrac{1+r}{1-r}\right)r\right)^m+\dfrac{2(1-\beta)r^N}{1-r}+4\mu sgn(t)\frac{ (1-\beta)^2r^Nt}{1-r}\\&\quad+\left(1+\frac{r}{1-r}\right)\frac{4\lambda (1-\beta)^2r^{2t+2}}{1-r^2}\leq \beta
	\end{align*}
	for $ |z|=r\leq R^{m, N, t}_{\mu, \lambda, \beta}(r) $, where $ R^{m, N, t}_{\mu, \lambda, \beta}(r) $ is smallest root in $ (0, 1) $ of $ \Phi^{m, N, t}_{\mu, \lambda, \beta}(r)=0 $, where $ \Phi^{m, N, t}_{\mu, \lambda, \beta} : [0, 1]\rightarrow\mathbb{R} $ defined by
	\begin{align*}
		\Phi^{m, N, t}_{\mu, \lambda, \beta}(r)&=\left(\beta r+(1-\beta)\left(\dfrac{1+r}{1-r}\right)r\right)^m+\dfrac{2(1-\beta)r^N}{1-r}+4\mu sgn(t)\frac{ (1-\beta)^2r^Nt}{1-r}\\&\quad+\left(1+\frac{r}{1-r}\right)\frac{4\lambda (1-\beta)^2r^{2t+2}}{1-r^2}-\beta.
	\end{align*}
	By a routine computation it can be easily shown that 
	\begin{align*}
		\frac{d\Phi^{m, N, t}_{\mu, \lambda, \beta}(r)}{dr}>0\; \mbox{for all}\; r\in (0, 1)
	\end{align*}
	which shows that the function $ \Phi^{m, N, t}_{\mu, \lambda, \beta} $ is increasing on $ (0, 1) $. Furthermore, we see that $ \Phi^{m, N, t}_{\mu, \lambda, \beta} $ is real-valued differential function on $ (0, 1) $ satisfying $ \Phi^{m, N, t}_{\mu, \lambda, \beta}(0)=-\beta<0 $ and $ \lim\limits_{r\rightarrow 1}\Phi^{m, N, t}_{\mu, \lambda, \beta}(r)=+\infty $, and this confirms, in view of intermediate value theorem, that $ R^{m, N, t}_{\mu, \lambda, \beta}(r) $ is the unique root in $ (0, 1) $ of $ \Phi^{m, N, t}_{\mu, \lambda, \beta}(r)=0 $. Certainly, we have 
	\begin{align}\label{B-eq-3.2}
		&\nonumber\left(\beta R^{m, N, t}_{\mu, \lambda, \beta}(r)+(1-\beta)\left(\frac{1+R^{m, N, t}_{\mu, \lambda, \beta}(r)}{1-R^{m, N, t}_{\mu, \lambda, \beta}(r)}\right)r\right)^m+4\mu sgn(t)\frac{(1-\beta)^2\left(R^{m, N, t}_{\mu, \lambda, \beta}(r)\right)^Nt}{1-R^{m, N, t}_{\mu, \lambda, \beta}(r)}\\&\quad+\frac{2(1-\beta)r^N}{1-R^{m, N, t}_{\mu, \lambda, \beta}(r)}+\left(1+\frac{R^{m, N, t}_{\mu, \lambda, \beta}(r)}{1-R^{m, N, t}_{\mu, \lambda, \beta}(r)}\right)\frac{ 4\lambda (1-\beta)^2\left(R^{m, N, t}_{\mu, \lambda, \beta}(r)\right)^{2t+2}}{1-\left(R^{m, N, t}_{\mu, \lambda, \beta}(r)\right)^2}=\beta.
	\end{align}
	For $f=h+\bar{g}\in \tilde{G}^0_{\mathcal{H}}(\beta)$, it can shown that 
	\begin{align}\label{Eq-3.8}
	\beta|z|+(1-\beta)\left(\frac{1-|z|}{1+|z|}\right)|z|	\leq |f(z)|\leq \beta|z|+(1-\beta)\left(\frac{1+|z|}{1-|z|}\right)|z|.
	\end{align}
	Therefore, in view of \eqref{Eq-3.8}, the Euclidean distance between $f(0)$ and the boundary of $f(\mathbb{D})$ is 
	\begin{align}\label{Eq-3.9}
		d(f(0), \partial f(\mathbb{D}))=\liminf_{r\to 1^{-}}|f(z)-f(0)|\geq \beta.
	\end{align}
	Therefore, in view of \eqref{Eq-3.9}, we see that the desired inequality $ S^f_{\mu, \lambda, m, N}(r)\leq \beta\leq d(f(0), \partial f(\mathbb{D})) $ is established for $ |z|=r\leq R^{m, N, t}_{\mu, \lambda, \beta}(r) $.\vspace{1.2mm}
	
	 The second step of the proof is to show that the constant $ R^{m, N, t}_{\mu, \lambda, \beta}(r) $ cannot be improved. Henceforth, we consider the function 
	\begin{align}\label{Eq-3.10}
		f_{\beta}(z)=z+\sum_{n=2}^{\infty}2(1-\beta)z^n.
	\end{align}
	Moreover, in view of \cite[Eq. (3.17)]{Allu-Hal-BSM-2021}, we see that 
	\begin{align}\label{eq-3.5}
		L(r)\leq |f_{\beta}(z)|\leq R(r)
	\end{align}
	and the equality holds also. From the left hand side of \eqref{eq-3.5}, a simple computation shows that $ d(f_{\beta}(0), \partial f_{\beta}(\mathbb{D}))=L(1)=\beta $. For $ r>R^{m, N, t}_{\mu, \lambda, \beta}(r) $ and $f=f_{\beta}$, we see that
	\begin{align*}
		&S^{f_{\beta}}_{\mu, \lambda, m, N}(r)\\&=\left(\beta r+(1-\beta)\left(\dfrac{1+r}{1-r}\right)r\right)^m+\dfrac{2(1-\beta)r^N}{1-r}+\frac{4\mu (1-\beta)^2r^{2t+2}}{1-r^2}\\&\quad+sgn(t)\left(1+\frac{r}{1-r}\right)\frac{4\lambda (1-\beta)^2r^Nt}{1-r}\\&>\nonumber\left(\beta R^{m, N, t}_{\mu, \lambda, \beta}(r)+(1-\beta)\left(\frac{1+R^{m, N, t}_{\mu, \lambda, \beta}(r)}{1-R^{m, N, t}_{\mu, \lambda, \beta}(r)}\right)r\right)^m+4\mu sgn(t)\frac{(1-\beta)^2\left(R^{m, N, t}_{\mu, \lambda, \beta}(r)\right)^Nt}{1-R^{m, N, t}_{\mu, \lambda, \beta}(r)}\\&\quad+\frac{2(1-\beta)r^N}{1-R^{m, N, t}_{\mu, \lambda, \beta}(r)}+\left(1+\frac{R^{m, N, t}_{\mu, \lambda, \beta}(r)}{1-R^{m, N, t}_{\mu, \lambda, \beta}(r)}\right)\frac{ 4\lambda (1-\beta)^2\left(R^{m, N, t}_{\mu, \lambda, \beta}(r)\right)^{2t+2}}{1-\left(R^{m, N, t}_{\mu, \lambda, \beta}(r)\right)^2}\\&=\beta\\&=d(f_{\beta}(0), \partial f_{\beta}(\mathbb{D}))
	\end{align*}
	which shows that $ R^{m, N, t}_{\mu, \lambda, \beta}(r) $ is best possible. This completes the proof. 
\end{proof}
\begin{proof}[\bf Proof of Theorem \ref{Th-3.6}]
	Let $ f=h+\bar{g}\in \tilde{G}^0_{\mathcal{H}}(\beta) $ for $ 0\leq\beta<1 $ be given by \eqref{B-eq-3.1}. An easy computation shows that 
	\begin{align*}
	r\varphi_0(r)+\sum_{n=2}^{\infty}(|a_n|+|b_n|)\varphi_n(r)&\leq r\varphi_0(r)+2(1-\beta)\sum_{n=2}^{\infty}\varphi_n(r)\leq \beta
	\end{align*}
	if $F_{\beta}(r)\leq 0$, where $F_{\beta}(r)$ is given by \eqref{Eq-3.4}. It is easy to see that $F_{\beta}(r)$ is a real valued differentiable function on $(0, 1)$. In view of \eqref{Eq-3.3}, we see that $F_{\beta}(0)<0$ and because $\{\varphi_n(r)\}$ is an increasing sequence of continuous function, we must have $\lim\limits_{r\to 1^{-}} F_{\beta}(r)=+\infty.$ Thus, it follows that $F_{\beta}(r)$ has a root in $(0, 1)$, say $R_{\beta}$. Moreover, $R_{\beta}$ is the unique root in $(0, 1)$ since 
	\begin{align*}
		R^{\prime}_{\beta}(r)=\varphi_1(r)+r\varphi_1^{\prime}(r)+2(1-\beta)\sum_{n=2}^{\infty}\varphi^{\prime}_n(r)>0\; \mbox{for all}\; r\in (0, 1).
	\end{align*}
	In view of \eqref{Eq-3.9}, we see that the desired inequality 
	\begin{align*}
		r\varphi_0(r)+\sum_{n=2}^{\infty}(|a_n|+|b_n|)\varphi_n(r)\leq d(f(0), \partial f(\mathbb{D}))
	\end{align*}
	holds for $r\leq R_{\beta}$. To show that $R_{\beta}$ is best possible, we consider the function $f_{\beta}$ given by \eqref{Eq-3.10}, and by the similar argument that being used in the proof of Theorem \ref{B-thm-3.1}, it can be shown using \eqref{Eq-3.4} for $r>R_{\beta}$ and $f=f_{\beta}$ that
	\begin{align*}
		r&\varphi_0(r)+\sum_{n=2}^{\infty}(|a_n|+|b_n|)\varphi_n(r)\\&>R_{\beta}\varphi_0(R_{\beta})+\sum_{n=2}^{\infty}(|a_n|+|b_n|)\varphi_n(R_{\beta})\\&=R_{\beta}\varphi_0(R_{\beta})+2(1-\beta)\sum_{n=2}^{\infty}\varphi_n(R_{\beta})\\&=F_{\beta}(R_{\beta})+\beta\\&=\beta\\&= d(f_{\beta}(0), \partial f_{\beta}(\mathbb{D}))
	\end{align*}
	which implies that $ R_{\beta} $ is best possible. This completes the proof.
\end{proof}	
	

\noindent\textbf{Compliance of Ethical Standards:}\\

\noindent\textbf{Conflict of interest} The authors declare that there is no conflict  of interest regarding the publication of this paper.\vspace{1.5mm}

\noindent\textbf{Data availability statement}  Data sharing not applicable to this article as no datasets were generated or analyzed during the current study.


\begin{thebibliography}{100}
	
	\bibitem{Aha-CMFT-2022} {\sc M. B. Ahamed}, The Bohr–Rogosinski Radius for a Certain Class of Close-to-Convex Harmonic Mappings, \textit{Comput. Methods Funct. Theory} (2022). https://doi.org/10.1007/s40315-022-00444-6.
	
	\bibitem{MBA-CVEE-2022} {\sc M. B. Ahamed}, The sharp refined Bohr–Rogosinski inequalities for certain classes of harmonic mappings, {\it Complex Var. Elliptic Equ.} (2022). https://doi.org/10.1080/17476933.2022.2155636.
	
	
	\bibitem{Aha-Aha-CMFT-2023} {\sc M. B. Ahamed} and {\sc S. Ahammed}, Bohr Type Inequalities for the Class of Self-Analytic Maps on the Unit Disk. Comput. Methods Funct. Theory (2023). https://doi.org/10.1007/s40315-023-00482-8.
	
	\bibitem{Aha-Allu-RMJ-2022} {\sc M. B. Ahamed} and {\sc V. Allu}, Bohr phenomenon for certain classes of harmonic mappings, \textit{Rocky Mountain J. Math.} 52(4): 1205-1225 (August 2022). DOI: 10.1216/rmj.2022.52.1205.
	
	\bibitem{Aha-Allu-CMB-2023} {\sc M. B. Ahamed} and {\sc V. Allu}, Bohr–Rogosinski radius for a certain class of close-to-convex harmonic mappings, \textit{Canad. Math. Bull.} (2023), 1-16.https://doi.org/10.4153/S0008439523000115
	
	\bibitem{Abu-CVEE-2010} {\sc Y. Abu-Muhanna},  Bohr's phenomenon in subordination and bounded harmonic classes, {\it Complex Var. Elliptic Equ.} {\bf  55} (2010), 1071--1078.
	
	\bibitem{Aizen-PAMS-2000} {\sc L. Aizenberg}, Multidimensional analogues of Bohr's theorem on power series, \textit{Proc. Amer. Math. Soc.} {\bf 128} (2000), 1147-1155.  
	
		\bibitem{Aizn-ST-2007} {\sc L. Aizenberg}, Generalization of results about the Bohr radius for power series, \textit{Studia Math.} {\bf 180} (2007), 161--168.
	
	\bibitem{aizenberg-2001} {\sc L. Aizenberg, A. Aytuna}  and {\sc P. Djakov}, Generalization of theorem on Bohr for bases in spaces of holomorphic functions of several complex variables, {\it J. Math. Anal. Appl.} {\bf  258} (2001), 429--447.

	
	
	\bibitem{Aleman-2012} {\sc A. Aleman} and {\sc A. Constantin}, Harmonic maps and ideal fluid flows, {\it Arch. Ration. Mech. Anal.} {\bf 204} (2012), 479--513.
	
	\bibitem{Alkhaleefah-Kayumov-Ponnusamy-PAMS-2019} {\sc S. A. Alkhaleefah, I. R. Kayumov}  and {\sc S. Ponnusamy}, On the Bohr inequality with a fixed zero coefficient, \textit{Proc. Amer. Math. Soc.} \textbf{147} (12) (2019), 5263-5274.
	
	\bibitem{Ali-Abdul-NG-CVEE-2016} {\sc R. M. Ali}, {\sc Z. Abdulhadi,} and {\sc Z.C. NG}, The Bohr radius for starlike logharmonic mappings,\textit{Complex Var. Elliptic Equ.} \textbf{61}(1), 1-14 (2016).
	
	\bibitem{Ponnusmy-Survey} {\sc R. M. Ali, Y. Abu-Muhanna}, and {\sc S. Ponnusamy, N.K. Govil}, et al. (Eds.), Progress in Approximation Theory and Applicable Complex Analysis, Springer Optimization and Its Applications, vol. 117 (2016), pp. 265-295.
	
	\bibitem{Allu-Arora-JMAA-2022} {\sc V. Allu} and {\sc V. Arora}, Bohr-Rogosinski type inequalities for concave univalent functions, \textit{J. Math. Anal. Appl.}(2022) 126845.
	
	\bibitem{Allu-Hal-BSM-2021} {\sc V. Allu} and {\sc H. Halder}, Bohr phenomenon for certain subclasses of harmonic mappings, \textit{Bull. Sci. Math.} \textbf{173}(2021), 103053.
	
	\bibitem{Allu-CMB-2022} {\sc V. Allu} and {\sc H. Halder}, Operator valued analogue of multidimensional Bohr inequality, {\it Canadian Math. Bull}. (2022).
	
	\bibitem{Allu-Hal-CMB-2022} {\sc V. Allu} and {\sc H. Halder}, Operator valued analogues of multidimensional Bohr’s inequality, \textit{Canad. Math. Bull.} \textbf{65}(4)(2022), 1020-1035.
	
	\bibitem{Allu-Hal-CMB-2023} {\sc V. Allu} and {\sc H. Halder}, Bohr operator on operator valued polyanalytic functions on simply connected domains, \textit{Canad. Math. Bull.} (2023), 1-11, https://doi.org/10.4153/S0008439523000541.
	
	\bibitem{Arora-CVEE-2022} V. Arora, Bohr’s phenomenon for holomorphic and harmonic functions with lacunary series in complex Banach spaces, \textit{Complex Var. Elliptic Equ.} (2022), https://doi.org/10.1080/17476933.2022.2146106.
	
	\bibitem{Aytuna-Djakov-BLMS-2013} {\sc  A. Aytuna} and {\sc P.  Djakov } Bohr property of bases in the space of entire functions and its generalizations, \textit{Bull. London Math. Soc.} \textbf{45}(2)(2013), 411-420.
	
	\bibitem{Beneteau-2004} {\sc C. B${\rm \acute{E}}$n${\rm \acute{E}}$teau, A. Dahlner} and {\sc D. Khavinson}, Remarks on the Bohr phenomenon, {\it  Comput. Methods Funct. Theory}  {\bf 4}  (2004), 1--19.
	
	\bibitem{Bhowmik-Das-PEMS-2021} {\sc B. Bhowmik} and {\sc N. Das}, , \textit{Proc. Edinburgh Math. Soc.} \textbf{64}(2021), 72-86.
	
	\bibitem{Boas-1997} {\sc H. P. Boas} and {\sc D. Khavinson}, Bohr's power series theorem in several variables, {\it Proc. Amer. Math. Soc.}  {\bf 125} (1997), 2975--2979.  
	
	\bibitem{Bohr-1914} {\sc H. Bohr}, A theorem concerning power series,  {\it Proc. Lond. Math. Soc.} s2-13 (1914), 1--5. 
	
	\bibitem{chen-Liu-Ponnusamy-arXiv-2023}  K. Chen, M. S. Liu, and S. Ponnusamy, Bohr-type inequalities for unimodular bounded analytic functions, \textit{Results Math.} \textbf{78}, 183 (2023).
	
	\bibitem{Constantin-2017} {\sc A. Constantin} and {\sc M. J. Martin}, A harmonic maps approach to fluid flows, {\it Math. Ann.} {\bf 369} (2017), 1--16.
	
	\bibitem{Das-JMAA-2022} {\sc N. Das}, Refinements of the Bohr and Rogosisnki phenomena, \textit{J. Math. Anal. Appl.} \textbf{508}(1) (2022), 125847.
	
	\bibitem{Das-CMB-2023} {\sc N. Das}, A logarithmic lower bound for the second Bohr radius, Canad. Math. Bull. (2023), 1-4. DOI: https://doi.org/10.4153/S0008439523000553
	
	\bibitem{Defant-Frerick-IJM-2001} {\sc A. Defant}, and {\sc L. Frerick}, A logarithmic lower bound for multi-dimensional Bohr radii, \textit{Israel J. Math} \textbf{152}(2006), 17-28.
	
	
	\bibitem{Dixon & BLMS & 1995} {\sc P. G. Dixon}, Banach algebras satisfying the non-unital von Neumann inequality, \textit{Bull. London Math. Soc.} \textbf{27}(4)(1995), 359--362.
	
	\bibitem{Evd-Pon-Ras-Antti-IM-2019} {\sc S. Evdoridis, S. Ponnusamy}, and {\sc A. Rasila}, Improved Bohr's inequality for locally univalent harmonic mappings, \textit{Indagationes Math.}\textbf{30}(2019) 201-213.
	
	\bibitem{Evdoridis-Ponn-RM-2021} {\sc S. Evdoridis, S. Ponnusamy}, and {\sc A. Rasila}, Improved Bohr’s inequality for shifted disks, \textit{Results Math.} 2021; 76(1):1–15.
	
	
	
	\bibitem{Galicer-Mansilla-Muro-TAMS-2020} {\sc D. Galicer}, {\sc M. Mansilla}, and {\sc S. Muro}, Mixed Bohr radius in several variables, \textit{Trans. Amer. Math. Soc.} \textbf{373} (2020),777-796.	
	
	\bibitem{Hamada-IJM-2009} {\sc H. Hamada}, {\sc T. Honda,} and {\sc G. Kohr}, Bohr's theorem for holomorphic mappings with values in homogeneous balls, \textit{Israel jour. Math} (2009) \textbf{173} 177-187.
	
	
	
	\bibitem{Kayumov-Khammatova-Ponnusamy-JMAA-2021} {\sc I. R. Kayumov, D. M. Khammatova}, and {\sc S. Ponnusamy}, Bohr-Rogosinski phenomenon for analytic functions and cesaro operators, {\it J. Math. Anal. Appl.} {\bf 496} (2021), 124824.
	
	\bibitem{Kayumov-Ponnusamy-JMAA-2018} {\sc I. R. Kayumov} and {\sc S. Ponnusamy}, Bohr’s inequalities for the analytic functions with lacunary series and	harmonic functions, \textit{J. Math. Anal. Appl.} \textbf{465}(2018), 857-871.
	
	\bibitem{Kayu-Kham-Ponn-MJM-2021} {\sc I. R Kayumov, D. M. Khammatova}, and {\sc S. Ponnusamy}, The Bohr inequality for the generalized Ces$ \acute{a} $ro	averaging operators, \textit{Mediterr J. Math.} 2021; 19: 15 pp.
	
	\bibitem{Kayumov-Ponnusamy-Shakirov-MN-2017} {\sc I. R. Kayumov}, {\sc S. Ponnusamy} and {\sc N. Shakirov}, Bohr radius for locally univalent harmonic mappings, \textit{Math. Nachri.} DOI: 10.1002/mana.201700068.
	
	\bibitem{Kumar-CVEE-2023} {\sc S. Kumar}, A generalization of the Bohr inequality and its applications, \textit{Complex Var. Elliptic Equ.} 68(6)(2023), 963-973.
	
	\bibitem{S. Kumar-PAMS-2022} {\sc S. Kumar}, On the multidimensional Bohr radius, \textit{Proc Amer Math Soc}, DOI: 10.1090/proc/16280.
	
	\bibitem{Kumar-JMAA-2023} {\sc S. Kumar} and {\sc R. Manna}, Revisit of multi-dimensional Bohr radius, \textit{J. Math. Anal. Appl.} \textbf{523}(2023), 127023.
	
	\bibitem{Kayumov-CRACAD-2018} {\sc I. R. Kayumov} and {\sc S. Ponnusamy}, Improved version of Bohr’s inequality, \textit{C. R. Acad. Sci. Paris, Ser.I} \textbf{356}(2018), 272--277.
	
	\bibitem{Lata-Singh-PAMS-2022} {\sc S. Lata} and {\sc D. Singh}, Bohr's inequality for non-commutative Hardy spaces, \textit{Proc. Amer. Math. Soc.} \textbf{150}(1) (2022), 201-211.
	
	\bibitem{Lewy-BAMS-1936} {\sc H. Lewy}, On the non-vanishing of the Jacobian in certain one-to-one mappings, \textit{Bull.
	Amer. Math. Soc.} \textbf{42} (1936), 689-692.
	
	\bibitem{Li-Ponn-CMJ-2016} {\sc L. Li} and {\sc S. Ponnusamy}, Injectivity of section of convex harmonic mappings and convolution theorem, \textit{Czechoslov. Math. J.} \textbf{66} (2016), 331-350.
	
	\bibitem{Lin-Liu-Ponnusamy-Acta-2023} {\sc R. Lin}, {\sc M. Liu},  and {\sc S. Ponnusamy}, The Bohr-type inequalities for holomorphic mappings with a lacunary series in several complex variables, \textit{Acta Math. Sci.} \textbf{43} (2023), 63-79.
	
	\bibitem{Gang-Liu-JMAA-2021} {\sc G. Liu}, Bohr-type inequality via proper combination, \textit{J. Math. Anal. Appl.} \textbf{503} (2021) 125308.
	
	
	\bibitem{Liu-Liu-JMAA-2020} {\sc X. Liu} and {\sc T-S Liu}, The Bohr inequality for holomorphic mappings with lacunary series in several complex variables, \textit{J. Math. Anal. Appl.} \textbf{485} (2020) 123844.
	
	\bibitem{Lin-Liu-Ponn-AMS-2023}{\sc  R. Lin, M. Liu}, and {\sc S. Ponnusamy}, The Bohr-type inequalities for holomorphic mappings with a lacunary series in several complex variables, \textit{Acta Math. Sci.}, \textbf{43B(1)}(2023), 63-79.
	
	\bibitem{Liu-Ponnusamy-BMMSS-2019} {\sc z. Liu} and {\sc S. Ponnusamy}, Bohr radius for subordination and $K$-quasiconformal harmonic mappings, \textit{Bull. Malays. Math. Sci. Soc.} \textbf{42},(2019), 2151–2168.
	
	\bibitem{Liu-Ponnusamy-PAMS-2021} {\sc M-S. Liu} and {\sc S. Ponnusamy}, Multidimensional analogues of refined Bohr's inequality, \textit{Proc. Amer. Math. Soc.} \textbf{149}(5), (2021), 2133-2146.
	
	
	\bibitem{Liu-Liu-Ponnusamy-2021} {\sc G. Liu}, {\sc Z. Liu,} and {\sc S. Ponnusamy}, Refined Bohr inequality for bounded analytic functions, \textit{Bull. Sci. Math.} \textbf{173} (2021), 103054.
	
	
	
	
	\bibitem{Paulsen-PLMS-2002} {\sc V. I. Paulsen}, {\sc G. Popescu,} and {\sc D. Singh},  On Bohr’s inequality \textit{Proc. Lond. Math.Soc.} \textbf{85}(2), 493-512 (2002).
	
	
	\bibitem{Ponnusamy-Vijayak-Wirths-RM-2020} {\sc S. Ponnusamy, R. Vijayakumar,} and {\sc K-J. Wirths}, New inequalities for the coefficients of unimodular bounded functions, \textit{Results Math} (2020) 75: Art 107
	
	
	\bibitem{Ponnusamy-Vij-Wirth-JMAA-2022} {\sc S. Ponnusamy, R. Vijayakumar}, and {\sc K.-J. Wirths}, Improved Bohr’s phenomenon in quasi-subordination classes, \textit{J. Math. Anal. Appl.} \textbf{506} (1) (2022), 125645.
	
	\bibitem{Rogosinski-1923} {\sc W. Rogosinski}, \"Uber Bildschranken bei Potenzreihen und ihren Abschnitten, \textit{Math. Z.}, \textbf{17} (1923), 260–276. 
	
	\bibitem{Wu-Wang-Long-RACSAM-2022} {\sc L. Wu, Q. Wang,} and {\sc B. Long,} Some Bohr-type inequalities with one parameter for bounded analytic functions, \textit{Rev. Real Acad. Cienc. Exactas Fis. Nat. Ser. A-Mat.} \textbf{116}, 61 (2022).
	
\end{thebibliography}
\end{document}